\documentclass[11pt]{amsart}
\usepackage{amsmath,amssymb,amsfonts}

\begin{document}

\title{\Large     On singular Sturm theorems  	}
\author{        D. Aharonov and U. Elias }
\date{} 

\begin{abstract}
The article summarizes some developments about a singular versions of the Sturm Comparison and 
Separation theorems where the coefficients or the interval of  definition may be unbounded.
\end{abstract}

\maketitle


The aim of this article is to summarize and popularize some results about 
Sturm comparison and separation theorems in singular cases.  

The classical Sturm's Comparison Theorem is formulated as following: 	\\
{\it  ``Consider the two differential equations 
\begin{align} 
		u'' & + \, p(x) \, u = 0,	\label{eq:Sturm-p}	\\
		v'' &  +   P(x)    v = 0,	\label{eq:Sturm-P} 
\end{align} 
where  $ P(x), p(x) $  are two continuous functions in an interval  $ [a,b] $, 
and  $ x_1, x_2 $  are two zeros of a solution  $ u $  of  (\ref{eq:Sturm-p})
there.  If 
$$ 
	P(x) \ge p(x)  \qquad  {\rm  but } \qquad  P(x) \not\equiv p(x) 
$$%
on  $ (x_1, x_2 ) $,   then every solution  $ v $  of  (\ref{eq:Sturm-P})  
has at least one  zero in  $ (x_1, x_2 ) $.''} 
  	
The Separation Theorem states that   {\it   ``if $ x_1, x_2 $  are two zeros of a 
solution  $ u $  of  (\ref{eq:Sturm-p}),  then every solution of the same equation 
which is linearly independent of  $ u(x) $  has a zero in  $ (x_1, x_2 ) $.''}

When the functions  $ P(x), p(x) $  are continuous only in an open interval but are 
unbounded at its endpoints and the zeros of the solution  $ u $  are located at 
the singular endpoints, the above formulations of Sturm's theorems are not valid.  
For example, consider 
$ 
	u'' + p_\lambda(x) u = 0 , 
$  
with  $  p_\lambda(x) = \dfrac{ 4\lambda(1 - \lambda) } { (1 - x^2)^2 } $, \  $ 0< \lambda < 1 $, 
which has the solutions  
$$ u_1(x) = (1-x)^\lambda (1+x)^{1 - \lambda} ,  \qquad	 u_2(x) = u_1(-x) $$%
that vanish at  $ x = 1, \, -1 $. 
If   $ 0 < \lambda < \mu < 1/2 $, then  $ p_\mu(x) > p_\lambda(x) $  on $ (-1,1) $, 
but the solution  $ \  v_1(x) = (1-x)^\mu (1+x)^{1 - \mu} \ $  of  $ \ v'' + p_\mu(x) v = 0 \ $ 
does not vanish between the zeros of  $ u_1(x) $.  Thus the claim of the Comparison Theorem 
is not satisfied.   The solutions  $ u_1(x), \ u_2(x) $   neither satisfy the Separation Theorem 
for the same equation.

The above difficulty may be answered by the following versions of the Sturm Theorems  which are valid 
for unbounded  $ p(x) $  and unbounded interval  $ (a,b) $:

\noindent
{\bf Singular Sturm Comparison Theorem. \ }{\it 
Consider the two differential equations  (\ref{eq:Sturm-p}), (\ref{eq:Sturm-P}) 
where  $ P(x), p(x) $  are two continuous functions  on the open, finite or infinite interval  
$ (a,b) $,   $ -\infty \le  a, \, b \le \infty $,   but not necessarily at its endpoints. 
Let the solution  $ u $  of  (\ref{eq:Sturm-p}) satisfy the boundary conditions 
\begin{equation}							\label{eq:principal}
	 	\int_{a} \frac{dx}{ u^2(x) } = \infty, 	  \qquad	
	   	\int^{b} \frac{dx}{ u^2(x) } = \infty .  	
\end{equation}	    
If    $ P(x) \ge p(x) $,  $ P(x) \not\equiv p(x) $  on  $ (a,b) $,  then every solution of the equation 
(\ref{eq:Sturm-P})    has a zero in  $ (a,b) $.
} 

\bigskip

\noindent
{\bf Singular Sturm Separation Theorem. \ }{\it 
Let  $ p(x) $  be continuous functions on  $ (a,b) $,   $ -\infty \le  a, \, b \le \infty $,   
but not necessarily at its endpoints. 
If  a solution   $ u $  of equation  (\ref{eq:Sturm-p})  satisfies the boundary conditions  
(\ref{eq:principal}),  then every solution of the same equation  which is linearly 
independent of  $ u(x) $  has a zero in  $ (a,b) $.	
} 

In the terminology of Hartman \cite[Chapter XI, Section 6]{Hartman},  boundary conditions  
(\ref{eq:principal})   mean that  $ u(x) $  is a  {\it principal solution}  of equation  (\ref{eq:Sturm-p}) 
at both endpoints of the interval,  i.e.,  $ u(x) $   is the essentially unique solution such that 
$  u(x) / \hat{u}(x) \to 0 $  as  $ x \to a $  and  $ x \to b $.

\noindent
{\bf Proofs.} \ The following argument is almost verbatim quoted from  \cite{AE}. 

Suppose that  $ u(x) > 0 $  in  $ (a,b) $,  otherwise apply the same arguments for a smaller interval. 
Take   $ c \in (a,b) $   so that   $ P(x) \not\equiv p(x) $  both in 
$ (a,c) $  and in  $ (c,b) $.  We show that the solution  $ v $  of  (\ref{eq:Sturm-P}) 
which is defined by the initial values  
\begin{equation}						\label{eq:ivp} 
		  v(c) = u(c) > 0, 	\qquad	v'(c)= u'(c) , 
\end{equation} 
has a zero in  $ (a,c) $  and a zero in  $ (c,b) $,  i.e.,  at least two zeros in $ (a,b) $.  
Suppose, on the contrary, that  $ v(x) \neq 0 $  in  $ (c,b) $,  and in fact, due to the 
initial conditions, $ v(x)>0 $  there.  From the identity  
	$ \  \left( v u' - u v' \right)' = v u'' - u v'' = (P-p) u v   \ $  
it follows that 
	$$  
		\left(v u' - u v' \right) \Big|_c^x  \  = \int_{c}^{x}  (P-p) u v  \, dx .	
	$$
By the initial value conditions  $ (v u' - u v')(c) = 0 $.  Since  $ u, v > 0 $   
and since  $ P \ge p, \ P(x) \not\equiv  p(x) $  in  $ [c,b) $,  the integral on the 
right-hand side increases and is positive for  $ c<x<b $.   So for some suitable value  $ d $,  
$ c < d < b $,   there exists a positive lower bound 
	$$ 
		(v u' - u v')(x) \ge C > 0,		\qquad		d \le x < b . 
	$$
Since  $ u \neq 0 $  in  $ [c,b) $, 
	$$  
		\left( -\frac{v}{u} \right)' = \frac{v u' - u v'}{u^2}  \ge  \frac{C}{u^2} > 0 , 
					       \qquad	 d \le x < b .  
	$$
Integration on  $ [d,x] $   yields 
\begin{equation}
	 	-\frac{ v(x) } { u(x) } + \frac{ v(d) } { u(d) } 
		  	\ge   \int_{d}^{x} \frac { C \, dx } { u^2(x) }   
\end{equation}
and by  (\ref{eq:principal}), 
	$$ \lim_{ x \to b^- }  \left[ -\frac{ v(x) } { u(x) }  \right] = +\infty ,  $$ 
contradicting the assumption that  $ u, v > 0 $  on  $ [c,b) $.   Thus  $ v $  must vanish
at least once in  $ (c,b) $.   By an analogous argument,  $ v(x) $  has another zero in  $ (a,c) $
as well.  Once  we know that the certain solution  $ v(x) $  of  (\ref{eq:Sturm-P}) 
has two zeros in  $ (a,b) $, it  follows by the classical separation theorem that every 
solution of  (\ref{eq:Sturm-P})  has at least one zero there.

The proof of the Separation Theorem is immediate: 
Without loss of generality we may assume that  $ u(x) \neq 0 $ in  $ (a,b) $,  otherwise we
may apply the same proof for a subinterval between two consecutive zeros.  Another linearly 
independent solution of  (\ref{eq:Sturm-p}) is   $  \  u(x) \int_{x_0}^x {dt} / { u^{2}(t) }  \  $   
and every solution is of the form 
\begin{equation}							\label{eq:gen-sol} 
	c_1 u(x) + c_2 u(x) \int_{x_0}^x \frac{dt}  { u^{2}(t) }
	= u(x)(	c_1 + c_2f(x) ) , 
\end{equation}  
with   $  \displaystyle   f(x) = \int_{x_0}^x \frac{dt} { u^{2}(t) } $.  \
According to  (\ref{eq:principal}),   
$ \displaystyle  \  \lim_{x \to a} f(x) = -\infty $,  
$ \displaystyle  \  \lim_{x \to b} f(x) =  \infty $, 
so every solution with  $ c_2 \neq 0 $  must change its sign  in  $ (a,b) $.

Only minor literal changes in the proof are needed to adjust the claims for the equations 
$ (ru')' + \, p(x) \, u = 0,	\  (rv')' +   P(x)    v = 0 $    with  $ r(x) > 0 $.  
\qquad$\square$


The boundary conditions   (\ref{eq:principal})  are also necessary for the validity of the 
Sturm Theorems in the following sense:

\noindent
{\bf Theorem.}{\sl \ 
Let   $ u(x) > 0 $  be a positive solution of equation   (\ref{eq:Sturm-p})  on $ (a,b) $,  
$ -\infty \le a < b \le \infty $,   and assume that at least one of the integrals 
\begin{equation}						\label{eq:L1L2}
   		L_1 := \int_a ^{x_0} \frac{dt}{u^2(t) } ,  \qquad 
   		L_2 := \int_{x_0} ^b \frac{dt}{u^2(t) } ,  \qquad a < x_0 < b, 
\end{equation}
%
  is finite.  	\\
(a) \ There exists  $ P(x) > p(x) $  such that the conclusion of Sturm Comparison Theorem 
does not hold for the pair of equations  (\ref{eq:Sturm-p}), (\ref{eq:Sturm-P}). 	\\
(b) \  The conclusion of Sturm Separation Theorem does not hold for equation  (\ref{eq:Sturm-p}) 
on the interval  $ (a,b) $. 
}

\noindent
{\bf Proof.} \ {\bf  The Comparison Theorem.} \     Let either one or both quantities 
$ L_1, \  L_2 $  be finite.  
Our aim is to construct a function  $ P(x) $  such that  $ P(x) > p(x) $  on  $ (a,b) $  
and some solution of the corresponding equation  (\ref{eq:Sturm-P})  has no zero in  $ (a,b) $.

 Let again  
$$ 
	f(x) = \int_{x_0} ^x \frac{dt}{u^2(t) }	
$$
and recall  that if  $ u(x) \neq 0 $  is a solution of  (\ref{eq:Sturm-p})  then 
$  \displaystyle  u(x) \int_{x_0}^x \frac{dt}{u^2(t) } = u(x) f(x)  $  is another solution of
the same equation.  Hence the function  $ f(x) $  is the ratio of two solutions.  
We shall utilize the Schwarzian derivative of  $ f(x) $,
$$  
    Sf = \left( {f''}/{f'} \right)' - \frac{1}{2} \left( {f''}/{f'} \right)^2 .
$$     
It is well known (and straightforward verified) that the Schwarzian of the ratio of any 
two solutions,  $ f = u_2/u_1 $,  satisfies  $  \  (Sf)(x) = 2 p(x) $.

Due to our assumption,  $ -L_1 \le f(x) \le L_2 $   and  $ (-L_1, L_2) $  may be 
a finite interval or a half ray.   We define
	$$ 	F(x) = g( f(x) )	$$%
where  $ g(t) $  is defined on   $ (-L_1, L_2) $,     
$ g'(t) > 0 $  and  $ Sg(t) > 0 $   there.  Suitable choices of   $ g(t) $  will be generated below.  
The Schwarzian derivative of this composition of mappings satisfies  
$$
	SF (x ) = S ( g \circ f ) (x ) = Sf (x) + \big( S g \circ f(x) \big) f'(x)^2 >  Sf(x) . 
$$
Now we define  $ P(x) $  by   $ SF(x) = 2P(x) $  and get accordingly a new equation  (\ref{eq:Sturm-P}) 
with   
$$
	P(x ) = \frac{1}{2} SF(x)  >  \frac{1}{2}Sf(x)  = p(x) . 
$$
Since  $ F'(x) = g'(f(x)) f'(x) > 0  \  $  and  $ F(x) $  increases, we may choose  $ v(x) > 0 $  so that 
$$    
	F(x) = \int_{x_0}^x \frac{dt}{v^2(t) } , 
$$%
and  $ v(x) = F'(x)^{-1/2} $  is a positive solution of  (\ref{eq:Sturm-P}) with no zero in
$ (a,b) $ .  Hence the Sturm Comparison Theorem does not apply to equations  (\ref{eq:Sturm-p}), (\ref{eq:Sturm-P}).

It remains only to generate a suitable  $ g(t) $  on  $ (-L_1, L_2) $, such that  $ g' > 0 $,  
$ Sg > 0 $.   
This is equivalent to choosing a differential equation  
$$
	w'' + R(t) w = 0   \qquad {\rm with}  \qquad   R(t) > 0 
$$  
which has a positive solution  $ w(t) $  on $ (-L_1, L_2) $ (i.e., it is disconjugate there). 
Indeed, for such $ w(t) > 0 $  we define  
	$ \displaystyle     g(t) = \int_{t_0} ^t w^{-2} $
and get  $ g'(t) = w^{-2}(t) > 0 $,  $  \  Sg(t) = 2 R(t) > 0 $.

If  $ L_1, L_2 $  are both finite, a simple choice is the differential equation 
$ 
	\displaystyle  w''(t) + \frac{ \pi^2 }{ (L_2 + L_1)^2 } \, w(t) = 0   \
$
with the solution  
$ \displaystyle   \  w(t) = \sin \left(  \frac{ \pi(t + L_1) }{ L_2 + L_1 } \right) > 0  \  $  
on   $ (-L_1, L_2) $.   Then 
$ \displaystyle   g(t) = -\cot \left(  \frac{ \pi(t + L_1) }{ L_2 + L_1 } \right) +$const,   \ 
$ \displaystyle   Sg(t) = \frac{ 2 \pi^2 }{ (L_2 + L_1)^2 } > 0 $,  and 
\begin{equation}						\label{eq:P>p-L1L2}
	P(x ) = p(x) + \frac{ \pi^2 }{ (L_2 + L_1)^2 } \, f'(x)^2  >  p(x) ,  
\end{equation} 
as required.

If  $ L_1 < \infty $,  $ L_2 = \infty $,   the equation 
$ 
	\displaystyle  w''(t) + \frac{ 1 }{ 4(t + L_1)^2 } \, w(t) = 0 $
with a solution  
	$ \displaystyle   w(t) = \sqrt{ t + L_1 } > 0 $  on   $ (-L_1, \infty) $ 
is suitable.  It leads to 
	$ g(t) = \log(t + L_1) + $const, \  $ Sg = { 1 }/{ 2(t + L_1)^2 } $, 
\begin{equation*} 
	P(x ) = p(x) + \frac{ 1 }{ 4(t + L_1)^2 } \Big|_{t=f(x)} f'(x)^2  >  p(x)   \qquad
	{\rm on}	\quad	(-L_1, \infty) . 
\end{equation*} 
This  $ g(t) $  is  applicable also for  $ L_2 < \infty $   if we restrict 
ourselves to  $ (-L_1, L_2) $.
The case of  $  \  -L_1 = -\infty $,   $ L_2 < \infty $  is treated similarly.

Note that if  $ -L_1 = -\infty $,  $ L_2 = \infty $,  then  no  equation  $ w'' + R(t) w = 0 $  
with  $ R(t) > 0 $  may have a positive solution on  $ (-\infty, \infty) $, 
since a concave function cannot be positive on the whole axis.

\noindent
{\bf The Separation Theorem.} \  
As mentioned in  (\ref{eq:gen-sol}),  the solutions of  (\ref{eq:Sturm-p})  are 
$$
   c_1 u(x) + c_2 u(x) \int_{x_0}^x \dfrac{dt} { u^{2}(t) }
	= c_1 u(x)( 1 + \frac{c_2}{c_1} f(x) )  
$$ 
and  $ -L_1 \le f(x) \le L_2 $. 
If   $ L_1 = -\infty $  and  $ L_2 < \infty $,  take  $  -1 / L_2 < c_2 / c_1 < 0 $;  
if  $ L_1 < \infty $,  $ L_2 = +\infty $,  take  $  0 < c_2 / c_1 < 1 / L_1 $; 
if both  $ L_1, L_2 < \infty $, take  $ c_2 / c_1 $  sufficiently small. 
In each case we get a solution of  (\ref{eq:Sturm-p})  which is linearly independent of  $ u(x) $  
and has no zero in  $ (a,b) $, contradicting the claim of the Separation Theorem.
\qed


\section*{Equivalent results}

An equivalent result had been proved earlier by Chuaqui et al in    \cite[Theorem 3]{Chuaqui-2009}, 
though  \cite{Chuaqui-2009}  studies different problems and applies different tools.   
The proof is based on the ``Relative Convexity Lemma''   
\cite[Lemma 2]{Chuaqui-2009}, \cite[p. 477] {Chuaqui-2007}  which may be stated as following: 
{\it 
Let  $ u $,  $v$  be respectively positive solutions of equations  (\ref{eq:Sturm-p})  and  (\ref{eq:Sturm-P}), 
$ P(x) \ge p(x) $.   Suppose that   
	$ \displaystyle   F(x) = \int_{x_0}^x \frac{dt}{ u^2 (t) } $.  
Then the function  $ w = \dfrac{v}{u} \circ{F^{-1}} $  is concave.%
}

In \cite[Theorem 3]{Chuaqui-2009}  it is proved that the boundary conditions  
(\ref{eq:principal})  are also necessary for the validity of the Sturm Comparison Theorem.  
The proof of the necessity in  \cite[p. 165]{Chuaqui-2009}  may be formulated as following.  
Multiplying   (\ref{eq:Sturm-p})  by  $ u $  and  (\ref{eq:Sturm-P})  by  $ v $ and substracting  
one from another:
\begin{equation}						\label{eq:v''u}
  -( P - p ) uv = v''u - vu'' = (v'u - vu')' = \left( u^2 \left(\frac{v}{u} \right)' \right)' . 
\end{equation}
With the choice  $ \displaystyle  \  P(x) = p(x) + \frac{ k^2 }{ u^4(x) } $,   \ 
equation  (\ref{eq:v''u})  becomes
\begin{equation}						\label{eq:v/u}
   u^2 \left( u^2 \left(\frac{v}{u} \right)' \right)' +k^2 \left(\frac{v}{u} \right) = 0 
\end{equation}
The change of variables   
	$ \displaystyle  \ y = f(x) = \int_{x_0}^x \frac{dt}{u^2(t) }   $  \
means 
$ \displaystyle  \  \frac{d}{dy} =u^2(x) \frac{d}{dx} $, \ thus equation  (\ref{eq:v/u})  is
$  \  \displaystyle  \frac{d^2}{dy^2} \left(\frac{v}{u} \right) +k^2 \left(\frac{v}{u} \right) = 0  \  $ 
and has a solution
\begin{equation*}
   \frac{v(x)}{u(x)} = \cos(ky) = \cos \left( k \int_{x_0} ^x \frac{dt}{u^2(t) } \right) .
\end{equation*} 
If both integrals in  (\ref{eq:principal})  are finite, then for sufficiently small  $k$, 
 $ v(x) \neq 0 $  on  $ (a,b) $.  Thus, if  (\ref{eq:principal}) is not satisfied,  
we found an equation   (\ref{eq:Sturm-P})  with    $ P(x) > p(x) $  for which the Comparison 
Theorem does not apply.

 Another equivalent conclusion is achieved implicitly also by Steinmetz in  \cite[Theorem 3]{Steinmetz}.  
\cite{Steinmetz}   studies the extension of univalent analytic functions on the unit disc,
therefore its results are formulated for a differential equation  $ y'' + q(x) y = 0 $  on 
$ [0, 1) $  and its symmetric extension to  $ (-1, 0] $.

The proof of necessity of (\ref{eq:principal}) in  \cite{Steinmetz}  is based on an argument 
of a different type.  
If  $ u_1, u_2 $  are two positive solutions of equation  (\ref{eq:Sturm-p}), 
then the function      $ v = u_1^\alpha u_2^{1 - \alpha} $, \  $ 0 < \alpha < 1 $,  
satisfies the differential equation 
\begin{equation}						\label{eq:u1u2}
   v'' + \left[ p(x) + \alpha( 1 - \alpha) \frac{ W^2(u_1, u_2) }{ u_1^2 u_2^2 } \right] v = 0 , 
\end{equation}
where  $ W(u_1, u_2) $  denotes the Wronskian,  which is in this case some constant number. 
((\ref{eq:u1u2})  is easily obtained if one utilizes the logarithmic derivative 
$ \displaystyle  \  \frac{v'}{v} = \alpha\frac{u_1'}{u_1} + (1 - \alpha)\frac{u_1'}{u_1}  \  $ 
and differentiate it once more).
If  $ u_1 $  is one positive solution of   (\ref{eq:Sturm-p}),  every other solution is of the form 
$ c_1 u_1(x) + c_2 u_1(x) \int_{x_0}^x u_1^{-2}(t) \, dt $.  So, if 
$ \int_{0}^1 {dt} / {u_1^2(t) } $  is finite,  a second positive solution 
$$
   u_2(x) = u_1(x) \left[ 1 - c \int_{0}^x \frac{dt}{u_1^2(t) }  \right] > 0  
$$
is chosen with  $ |c| $  is sufficiently small.  Utilizing the two positive solutions  $ u_1(x), u_2(x) $, 
\cite{Steinmetz}  obtains according to  (\ref{eq:u1u2})  an equation  (\ref{eq:Sturm-P})  with 
\begin{equation}						\label{eq:P-u1u2} 
    	P(x) = p(x) + \frac{ k^2 }{ u_1^2 u_2^2 }  >  p(x)
\end{equation}
and a positive solution  $ v = u_1^{\alpha} u_2^{1 - \alpha} $   on  $ (a,b) $   for which the claim 
of the Comparison Theorem does not apply.




\end{document}